\pdfoutput=1

\documentclass[11pt]{article}

\usepackage[final]{acl}

\usepackage{times}
\usepackage{latexsym}

\usepackage[T1]{fontenc}

\usepackage[utf8]{inputenc}

\usepackage{color}
\usepackage{algorithm}
\usepackage{algorithmic}
\usepackage{amssymb}

\usepackage{booktabs}
\def\toprule{\noalign{\hrule height 1.2pt\smallskip}}
\def\midrule{\noalign{\smallskip\hrule\smallskip}}
\def\bottomrule{\noalign{\smallskip\hrule height 1.2pt}}

\usepackage{pgfplots}
    \definecolor{turbo1}{rgb}{0.298 0.573 0.980} 
    \definecolor{turbo2}{rgb}{0.882 0.290  0.125} 

\usepackage{microtype}

\usepackage{inconsolata}

\title{Solving General Natural-Language-Description Optimization Problems with Large Language Models}

\author{Jihai Zhang, Wei Wang, Siyan Guo, Li Wang, Fangquan Lin, Cheng Yang \and Wotao Yin \\
        Alibaba Group \\ \texttt{\{jihai.zjh, zhuazhua.ww, guosiyan.gsy, feiyu.wl, } \\ \texttt{fangquan.linfq, charis.yangc, wotao.yin\}@alibaba-inc.com}}

\begin{document}
\maketitle
\begin{abstract}
Optimization problems seek to find the best solution to an objective under a set of constraints, and have been widely investigated in real-world applications. Modeling and solving optimization problems in a specific domain typically require a combination of domain knowledge, mathematical skills, and programming ability, making it difficult for general users and even domain professionals. 
In this paper, we propose a novel framework called \textit{OptLLM} that augments LLMs with external solvers.
Specifically, OptLLM accepts user queries in natural language, convert them into mathematical formulations and programming codes, and calls the solvers to calculate the results for decision-making. In addition, OptLLM supports multi-round dialogues to gradually refine the modeling and solving of optimization problems. 
To illustrate the effectiveness of OptLLM, we provide tutorials on three typical optimization applications and conduct experiments on both prompt-based GPT models and a fine-tuned Qwen model using a large-scale self-developed optimization dataset. 
Experimental results show that OptLLM works with various LLMs, and the fine-tuned model achieves an accuracy boost compared to the prompt-based models.
Some features of OptLLM framework have been available for trial since June 2023~(https://opt.alibabacloud.com/chat or https://opt.aliyun.com/chat).
\end{abstract}

\section{Introduction}
Optimization problems have been widely investigated in real-world domains including financial investment~\cite{ye2020financialinvestment}, supply chain management~\cite{li2023supplychain}, logistics transportation~\cite{xie2020logistictransport} and competitive strategy~\cite{silver2017competitivestrategy}.
Such ubiquitous optimization problems raise critical demands for efficient modeling and solving methods.

Currently, modeling and solving optimization problems in a specific domain usually involves three steps~\cite{ramamonjison2022optimization}.
First, based on domain knowledge, experts summarize the application scenarios into problem descriptions using natural language or mathematical formulas, with clear indication of variables, objectives, constraints, and parameters. Second, experts extract and encode critical information from the problem descriptions with modeling languages such as Python, R or AMPL. Finally, the optimization process is carried out by experts or solvers to obtain the final decision-making results. 
Meanwhile, the entire process calls for a combination of domain knowledge, mathematical skills, and programming ability, which is unfriendly to beginners or even professionals in that domain.

Recently, large language models (LLMs) have demonstrated strong capabilities in natural language understanding and generation ~\cite{openai2023gpt4}. 
However, despite LLMs' strong performance across a range of NLP tasks (\textit{e.g.}, content generation and Q\&A dialogue)~\cite{brown2020gpt3}, their ability in arithmetic and logical reasoning may be insufficient and unfaithful~\cite{imani2023mathprompter}. On the other hand, data privacy remains one concern for online services like GPT-4~\cite{openai2023gpt4}. That is, inclusion of domain-specific information in prompts may cause data breach at the LLM service provider side or during transmission in public networks, even under the service level agreements for privacy~\cite{li2023supplychain}. Hence, deployment of open-resourced LLMs (\textit{e.g.}, Llama~\cite{touvron2023llama}, PaLM~\cite{anil2023palm}, and Qwen\footnote{https://modelscope.cn/models/qwen/Qwen-7B-Chat}) is preferred for privacy-sensitive applications.

In light of these above, we propose \textit{OptLLM}, a framework unifying either open-sourced LLMs or online LLM services, and external solvers for automated modeling and solving of optimization problems.
Specifically, OptLLM consists of three modules. 
First, the interaction refinement module interacts with users to complete problem descriptions and ensures the input is a valid optimization problem.
Second, the converter module converts problem descriptions to mathematical formulations and programming codes, and ensures the codes are correct.
Last, the responser module sends the code to an external optimization solver, receives the results and interprets them. OptLLM allows users to iteratively refine any stage outputs through chatting or direct editing, until satisfactory results are obtained. In this way, OptLLM aims to make it significantly easier for users to model and solve optimization problems.


\section{Related Work}
LLMs, or large language models, are predominantly Transformers~\cite{NIPS2017transformer} trained on extensive text corpus from various sources (\textit{e.g.}, webs and books~\cite{brown2020gpt3}). They are trained to predict the next token in a given context, and could generate coherent responses after fine-tuning and alignment~\cite{openai2023gpt4}. Below we briefly cover applications and techniques related to automated optimization problem modeling and solving using LLMs.

\subsection{Applications of LLMs}
With the widespread attention on LLMs, their applications are popping up in varied domains, such as open-domain Q\&A~\cite{liu2023prompt}, database management~\cite{zhou2023llm}, and strategizing agents~\cite{yao2022react}. Studies on arithmetic reasoning, or mathematical reasoning~\cite{qiao2022reasoning}, investigate the ability of LLMs to solve math word problems (MWP)~\cite{patel2021mwp}.
Existing work mainly focus on general math problems including function evaluation, numerical calculation and theorem 
proving~\cite{imani2023mathprompter,yang2023leandojo}. Unfortunately, the reasoning ability of LLMs is still far from being usable~\cite{wang2022selfconsistency} and even competent models like GPT-4 are inconsistently bad at numeric calculations. In contrast, our work relies on LLMs to model optimization problems, and external solvers for solving them. 

    
    
    

\subsection{Techniques of LLMs}
To adapt LLMs to downstream tasks, two strategies are commonly used: \textit{prompting} and \textit{supervised fine-tuning} (SFT)~\cite{liu2023prompt}.
Prompting, also known as in-context learning, leverages additional task information, zero to a few domain-specific examples, and expected answer format to guide LLMs without additional training. 
Recent works show that specially-designed prompts, such as those via chain-of-thoughts~\cite{wei2022chain}, iterative refinement~\cite{madaan2023selfrefine} and black-box prompt tuning~\cite{sun2022black}, can significantly improve the performance of LLMs on downstream tasks. 
On the other hand, SFT leverages task-specific data and objective functions to train LLMs, which demonstrates a significant enhancement in downstream applications~\cite{baldazzi2023finetuning}. SFT is more effective than prompting when such task-specific data are available. 

\begin{figure*}[htbp]
\centering
\includegraphics[width=0.9\textwidth]{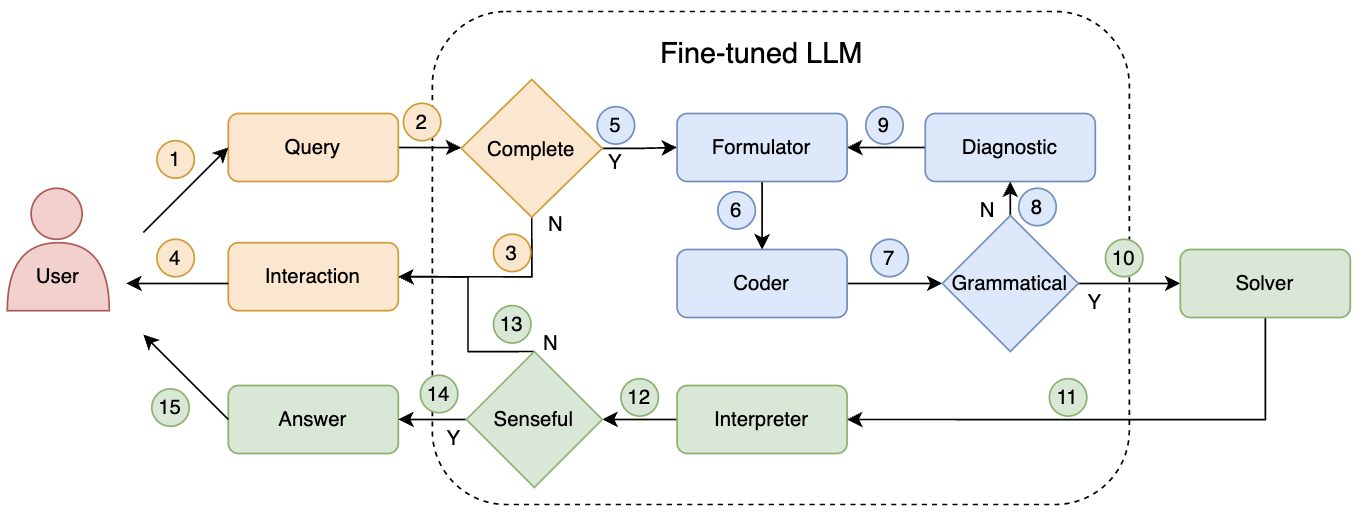} 
\caption{OptLLM framework consists of three main modules: (1) \textit{Interaction Refinement Module}, Step 1 to 4 (marked in orange), interact with the user to get a complete problem description in natural language; (2) \textit{Converter Module}, Step 5 to 9 (marked in blue), converts problem description to math formulas and codes; and (3) \textit{Responser Module}, Step 10 to 15 (marked in green), calls the solver, checks and interprets its results, and responses to the user.}
\label{framework}
\end{figure*}

\section{Proposed Framework: OptLLM}

We propose OptLLM that unifies LLMs and external solvers for automated modeling and solving of optimization problems. By designing OptLLM to interact with domain users via natural language, we hope to reduce the need for specialized knowledge on optimization or coding, and improve the experience for end-users. 
OptLLM primarily consists of three modules: interaction refinement module, converter module, and responser module.

\subsection{Interaction Refinement Module}
As shown in Figure~\ref{framework}, the interaction refinement module consists of Step $1$ to $4$ (marked in orange). Step $1$, the user queries OptLLM in natural language. Step $2$, the queries are pre-processed, including inserting instructions and prompt engineering. The pre-processing is used to clarify the task and output formats for LLMs. For online LLM service like GPT, a typical instruction could be ``You are an operation research expert and your task is to model the optimization problem given its description in natural language.'' The queries are then checked in the `Complete' part. Complete queries should have clear indication of variables, objectives, constraints, and parameters for optimization. If the user's queries are complete, the queries are sent to the next module for modeling. Otherwise, OptLLM detects some information is missing and request user to provide more details. We will provide an example application in Application 2 below. In practice, OptLLM responds to user inputs in various scenarios. If a user's queries are unrelated to optimization problems, OptLLM would prompt and guide the user towards asking optimization related questions. 

\subsection{Converter Module}
The Converter module contains Step $5$ to $9$ in Figure~\ref{framework}  (marked in blue). The module is used to convert problem descriptions in natural language to codes and check their grammar. Step $5$ receives the output of the interaction refinement module and passes it to the `Formulator' of OptLLM. The Formulator translates natural language descriptions into the corresponding formulas for objective and constraints. Then in Step $6$, the formulas are fed into `Coder' to generate corresponding code in a preset programming language. In `Grammatical', the code will be checked for grammar mistakes. If the syntax test fails, it enters the diagnostic module and OptLLM reformulates it based on its own feedback. Otherwise, the code will be sent to an external solver.

We use MindOpt Algebraic Programming Language, or MAPL\footnote{https://www.yuque.com/mindopt/apl\_en/tuhebr} as the default programming language. Designed by Alibaba, MAPL is an efficient and versatile modeling language that supports many mainstream solvers, including MindOpt, Gurobi, CPLEX, Ipopt, Cbc. We use MindOpt\footnote{https://opt.aliyun.com/} by default.

\begin{figure*}[h]
\centering
\includegraphics[width=0.80\textwidth]{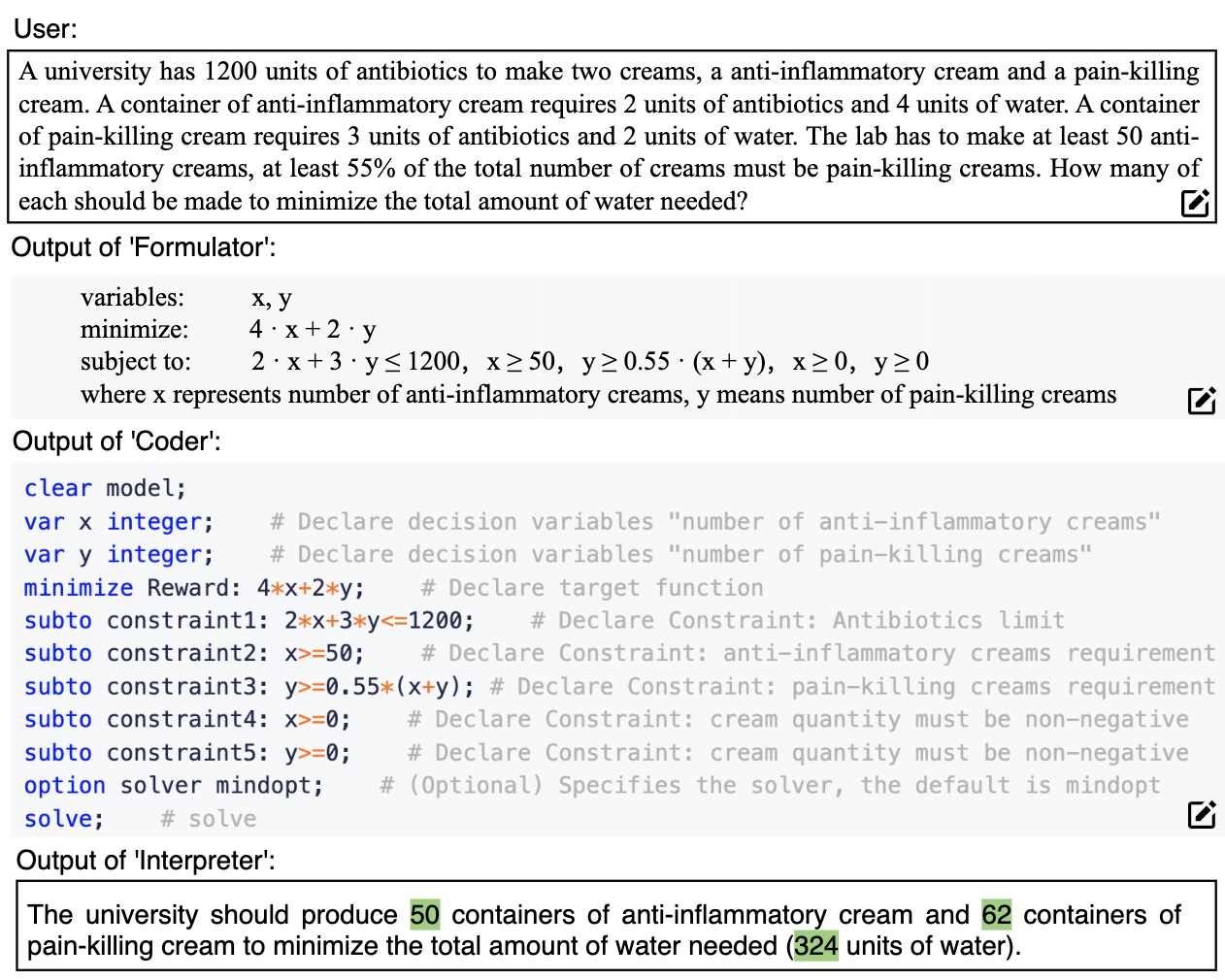} 
\caption{Overview of Application 1. The user provides a complete description for an optimization problem. The `Formulator' translates it into formulas, then `Coder' generates the corresponding MAPL code. At last, `Interpreter' receives the solver output and interprets it with natural language. The user input, formulas and code can be directly edited and the rest parts will be re-generated.}
\label{app1}
\end{figure*}

\subsection{Responser Module}
The Responser Module consists of Step $10$ to $15$ in Figure~\ref{framework}  (marked in green). In Step $10$, the programming code is sent to the solver. In Step $11$, the `Interpreter' block collects the solver's solution and interprets it in natural language. The solution and interpretation are then subjected to `Senseful', which checks semantic validity. A solution fails if it does not meet a user-defined requirement (e.g. the user requires an integer solution but the solver returns a real value). 
If a solution fails, the issue is resolved through interaction with the user. Otherwise, the solution and interpretation will be formatted by the `Answer' block and presented to the user.


\section{Applications}
Our framework can solve generic optimization problems based on their natural language description. In this part, we introduce three basic applications, including single-round QA with complete description, multi-round conversations with missing information detection, and optimization problem solving with external data.

\subsection{Application 1: Single-round QA}
In the single-round QA application, we assume the user has provided a complete natural language description of the optimization problem such that the variables, objective and constraints can be deduced. This application is often used in the education, e.g., when a student enters a complete optimization problem into the system, or when a teacher lectures a student with a complete problem. 
As show in Figure~\ref{app1}, the `Formulator' block generates the corresponding formulas, with the ability to automatically detect variable names which are not explicitly specified in the problem description; the `Coder' then generates the corresponding MAPL codes based on the contextual formulas and calls the `Solver' for solving. Finally, the `Interpreter' translates the solver's solution into natural language, making it easy for the user to understand. Each module mentioned above can be optionally displayed or manually edited by the user. For example, if the user wants to make a numerical substitution or change one constraint, it can be quickly achieved through the editing function. Once modified, OptLLM would re-generate the rest parts accordingly.

\subsection{Application 2: Multi-round Conversations}
\label{sec:application_2}
The Application 1 assumes user has provided a complete problem description. In many scenarios, users may not provide such a description at once, especially if they would like to gradually build up a complex problem. In light of this, it is necessary to guide them step by step through interactions to provide the necessary information for modelling optimization problems. We hope to start with the simplest chat, detect the missing information, and gradually guide users through interaction to provide necessary information indicating the variables, objective and constraints. 

\begin{figure}[h]
\centering
\includegraphics[width=0.48\textwidth]{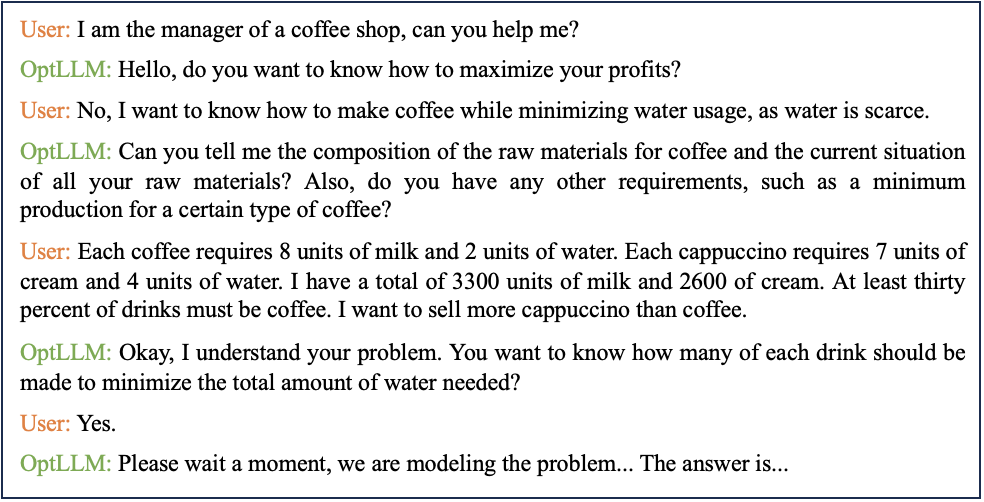} 
\caption{Overview of Application 2. OptLLM attempts to guide user to provide the necessary information for an optimization problem and then provide the answer directly to the user. The math formulas and codes are hidden.}
\label{app2}
\end{figure}

Figure~\ref{app2} shows an example where a coffee shop owner uses the system without any knowledge of optimization. The model guides the user to provide the objective and necessary constraints, and finally provides the answer. For people without a math or coding background, we provide options to hide the formulas, codes and intermediate processes. The system enables the users to enjoy the benefits of programming language and solvers within natural language dialogues.

\subsection{Application 3: External Data Files}
There are scenarios where the data for an optimization problem cannot be concisely tabulated or embedded in the problem description. In addition, LLMs typically have a token limit of a few thousands, which could be easily exceeded by the lengthy descriptions or multiple rounds of interactions, if lots of data are embedded. To address this, we design OptLLM to accept external data files with a predetermined format. Users may query the system with instructions on from which files each part of the data can be acquired. Inspired by LangChain\footnote{https://github.com/langchain-ai/langchain}, the data will not be passed to the model to save tokens and further preserve data privacy. Instead, only the external solver will access the data files in order to calculate the final solutions. 

\section{Fine-tuning Large Language Model}
OptLLM permits both online LLM service or open-sourced LLMs as the base model. In this section, we introduce the fine-tuning process on Alibaba's self-developed Qwen model. Considering the large size of the Qwen model we use (50B parameter version) and a limited budget (eight NVIDIA V100 GPUs), the data scale and existing hardware do not support continuous pre-training or full parameter fine-tuning. Thus, we adopt LoRA (low-rank adaption)~\cite{hu2021lora}, a parameter-efficient fine-tuning scheme (PEFT)~\cite{houlsby2019parameter}.

\begin{figure}[h]
\centering
\includegraphics[width=0.45\textwidth]{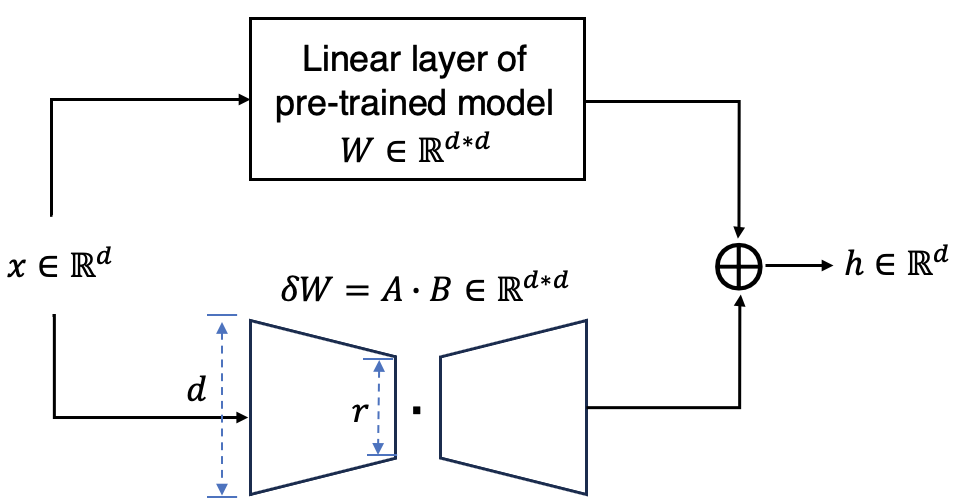} 
\caption{Overview of LoRA (low-rank adaption).}
\label{ft}
\end{figure}
As shown in Figure~\ref{ft}, at each linear layer of the Qwen model, LoRA inserts two trainable low-rank matrices $A\in \mathbb{R}^{d\times r}$ and $B\in \mathbb{R}^{r\times d}$ to approximately optimize the original parameters: $W^{new} = W + A\cdot B$, where $\textit{x}$ is our fine-tuning data, and $\textit{h}$ is the output of the linear layer, $W$ is the fixed original parameter matrix. Overall, the amount of parameters introduced by LoRA is below $1\%$ of the original model.

\section{Experiments}

\subsection{Datasets}
We focus primarily on linear programming (LP) and mixed integer linear programming (MILP) problems, which may be of strong interests in industrial applications. To the best of our knowledge, there is currently no publicly available datasets on general optimization problems except the data from NL4OPT competition\footnote{https://nl4opt.github.io/}~\cite{ramamonjison2023nl4opt}. Few additional problems can be crawled from websites, but they may still not be sufficient for fine-tuning the LLMs. Thus, we constructed our own fine-tuning and test datasets. We ensure that Qwen model has not seen the test datasets during the pre-training phase.

\begin{figure}[h]
\centering
\includegraphics[width=0.48\textwidth]{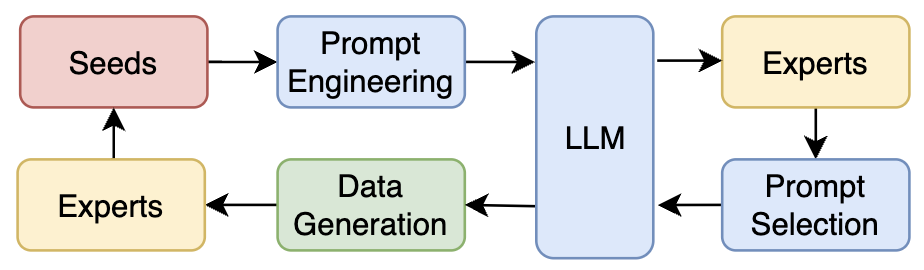} 
\caption{Flowchart for data collection.}
\label{traindata}
\end{figure}

\textbf{Fine-tuning Dataset.} Figure~\ref{traindata} shows the data collection process. To build a large-scale dataset, we start with seed optimization problems manually designed by experts, followed by designing prompts and calling LLMs to generate more problems. We then manually label the data and select prompts that perform well. The resultant prompts are used to generate more data, which are again manually labelled. The labeled data can be used as seeds for the next round of data generation. This process is repeated several times. Finally, we collected a high-quality optimization training dataset with total 15k instances in English and Chinese.

\textbf{Test Dataset.} We select $100$ optimization problems with natural language description in English, called En$100$, part of which are from the dev dataset of NL4OPT competition. 
Besides, we also prepare another $100$ Chinese optimization problems, called CN$100$. All test data are manually checked to ensure correctness.

\subsection{Metrics}
In this study, we focus on evaluating the model performance on single-round QA as in Application 1. The multi-round conversations could be significantly more diverse than single-round QA, and may have multiple acceptable answers. We wish to evaluate it fairly and faithfully in the future when our OptLLM is fully deployed. As for single-round QA, it should be noted that different formulas can lead to the same solution, e.g., adding a redundant constraint $x\ge10$ to an existing one $x\ge20$ will not affect the solution but $x\ge10$ may have not been mentioned by the problem description. Thus, evaluating the model by the solver results may overlook the mistakes in formulas, and we propose to evaluate the model by accuracy in formula generation. The model is considered correct on one sample only when all the variables, objective and constraints it generates matches exactly with the ground truth on that sample. 
It should be noted that our metric is more strict than the declaration-level mapping accuracy used in NL4OPT study~\cite{ramamonjison2023nl4opt}. 

\subsection{Implementation Details}
We compare the finetuned Qwen model against two prompt-based models: GPT-3.5 (gpt-3.5-turbo)~\cite{ouyang2022training} and GPT-4~\cite{openai2023gpt4} under our OptLLM framework. For GPT-3.5 and 4, we use the standard one-shot prompt: ``You are an expert in mathematical programming. Please refer to Case 1 and provide a JSON expression for Problem 1 with explanations. Case1: \{Question\_and\_Answer\_of\_Case1\}, Problem1: \{Question\}." We have also tried prompts with more shots but the performance does not improve significantly, so we stick to one shot. For Qwen model, we finetune it using LoRA as describe in previous section. LoRA is inserted at every linear layer of the model. The dimension $r$ of all LoRA layers is set to $32$. The AdamW optimizer is used with an initial learning rate of $0.0002$, $\beta=[0.9,0.999]$, and a linear decay schedule. The number of training epochs is set to $20$ with a mini-batch size of $32$ due to limited GPU memory. The model is implemented under HuggingFace's Transformers library~\cite{shen2023hugginggpt} and trained on eight NVIDIA V100 GPUs using DeepSpeed Zero stage 3~\cite{yao2023deepspeed}.

\subsection{Results}

\begin{table}[h]
\centering
\caption{The accuracy of LLMs on test datasets.}
\begin{tabular}{lccc}
\toprule
Datasets & GPT-3.5 & GPT-4 & Qwen-SFT \\
\midrule
EN100    & 71\%      & 82\%    & 87\%     \\
CN100    & 71\%      & 80\%    & 80\%    \\
\bottomrule
\end{tabular}
\label{table1}
\end{table}

\textbf{Overall performance.} 
As shown in Table~\ref{table1}, the supervised fine-tuning Qwen, or Qwen-SFT, surpassed GPT-3.5 on both datasets. It also achieved comparable performance to GPT-4 on CN$100$ and exceeded GPT-4 on EN$100$. We manually identified the specific error causes - GPT-3.5 and GPT-4 made mistakes in identifying strict constraints, e.g., the problem description states ``A is more than B", that is, ``$A > B$", but both GPTs inferred ``$A \geq B$". In contrast, Qwen-SFT had more successes in identifying such constraints, owing to the fine-tuning process enabling it to learn sophisticated patterns. 



\begin{figure}[h]
    \centering
    \begin{tikzpicture}
    \begin{axis}[
        xlabel={Fine-tuning Epochs},
        ylabel={Test Accuracy (\%)},
        xmin=0, xmax=20,
        ymin=50, ymax=100,
        xtick={0,4,8,12,16,20},
        ytick={50,60,70,80,90,100},
        legend style={draw=none},
        legend pos=south east,
        ymajorgrids=true,
        grid style=dashed,
        width=0.9\linewidth, height=0.54\linewidth
    ]
    
    \addplot[
        color=turbo1,
        line width=1.6pt,
        mark=square,
        ]
        coordinates {
        (1,56)(2,72)(4,74)(6,83)(8,80)(10,87)(12,86)(14,85)(16,88)(18,90)(20,87)
        };
    \addplot[
        color=turbo2,
        line width=1.6pt,
        mark=o,
        ]
        coordinates {
        (1,55)(2,59)(4,62)(6,80)(8,73)(10,78)(12,75)(14,78)(16,77)(18,82)(20,80)
        };
    \legend{EN100,CN100}
        
    \end{axis}
    \end{tikzpicture}
    \caption{Test Accuracy at different finetuning epochs.}
    \label{ab1}
\end{figure}
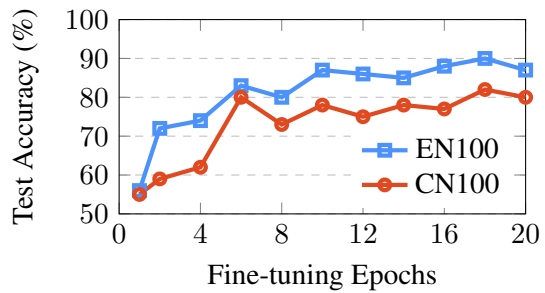

\textbf{Impact of fine-tuning epochs.} Figure~\ref{ab1} shows that, the model's performance on the test datasets improves with more fine-tuning epochs, and start to plateaus after 10 epochs. Given the prolonged training time, we set the fine-tuning epochs to be 20 by default.

\begin{table}[h]
\centering
\caption{The impact of finetuning data size.}
\begin{tabular}{lcccc}
\toprule
\#Samples    & 500      & 1000    & 2000  & 4000    \\
\#Epoch    & 40      & 20    & 10  & 5   \\
\midrule
Accuracy    & 28\%   & 42\%    & 50\% & 57\%    \\
\bottomrule
\end{tabular}
\label{ab2}
\end{table}

\textbf{Impact of data diversity.} To investigate the influence of fine-tuning data size on model performance, we vary the number of samples and epochs so that, in each setting, the model is trained on roughly the same number of tokens. We fine-tune and evaluate the model on Chinese data only. As shown in Table~\ref{ab2}, the model performance increases along with the data size. This indicates that we should collect as many diverse data as possible to achieve better results.

\section{Path to Deployment}

\begin{figure}[h]
\centering
\includegraphics[width=0.45\textwidth]{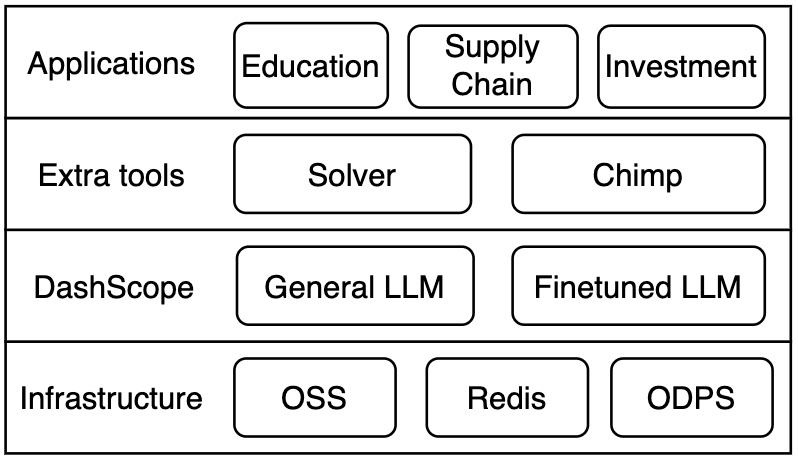} 
\caption{Overview of the deployment framework.}
\label{deploy}
\end{figure}

The proposed OptLLM framework can be deployed on the cloud. We take Alibaba Cloud\footnote{https://www.aliyun.com/} as an example to illustrate the deployment of OptLLM. As shown in Figure~\ref{deploy}, the infrastructure includes: i) the OSS provides data storage for user data that may be used in Application 3; ii) Redis is used for recording online conversation context; iii) ODPS is used for logging historical logs.

DashScope\footnote{https://dashscope.aliyun.com/} is an inference platform that supports both existing LLM APIs~(\textit{e.g.}, GPT-3.5 and GPT-4) or self-built LLMs~(\textit{e.g.}, Qwn-SFT and Llama2-SFT). External tools include Solver, such as MindOpt, and Chimp, a testing platform for the entire framework. Once deployed on the cloud, the proposed OptLLM framework has the potential to support applications in various domains, such as educational services, financial investment and supply chain management. In June 2023, we have deployed the first version on Alibaba Cloud, which includes some of the features introduced in this paper, with more features currently under development.

\section{Limitation}
Although our system is capable of handling single-round optimization problems, as well as multi-round addition, deletion, and modification operations for some optimization problems, our model's effectiveness will be somewhat affected when dealing with incomplete issues that require additional knowledge for certain parts. This is because this extra knowledge may not be possessed by our large model due to certain reasons, such as our model's knowledge base being up-to-date only until 2023, meaning it wouldn't be aware of knowledge from 2024. There are two ways to address this issue: one is to update the underlying large model in real-time, but this would entail significant financial and material costs. The other option involves using methods related to Retriever-Augmented Generation (RAG). These are aspects we plan to explore in our future work.

\section{Conclusion}
In this paper, we propose OptLLM, an effective framework that augments LLMs (such as Qwen model and GPT-4) with external solvers for automated modeling and solving of optimization problems. Specifically, OptLLM comprises three modules:
the interaction module for completing the problem description, the converter for translating the description into code, and the responser for calling solvers and interpreting the results, respectively. By iterating the above steps through chatting with users, OptLLM has the potential to assist both beginners and domain professionals to achieve faithful decision-making for optimization problems. We illustrate the effectiveness of OptLLM with three proof-of-concept applications and experiments.
In the future, we will focus on promoting the diversity of optimization problems by including more real-world cases from various domains and scenarios. We will also explore methods to enhance arithmetic and logical reasoning, as well as more open-sourced LLMs and evaluation methods.




\section{Acknowledgments}
We would like to express our sincerest gratitude to the anonymous reviewers for their insightful feedback on our work. We are also immensely thankful to our colleagues: Hu Jiang, You Wu, Churan Liu, Binyang Shen, Junqiu Pan, Mou Sun, Jiwei Li, Ao Zhang, Yuhua Song, Liang Zhao, Wei Jiang, Zhongkai Yi and Hanwei Zhang for their invaluable support throughout the research process. We would also like to extend our heartfelt thanks to Professor Zaiwen Wen's team and Professor Zhifang Yang's team for their assistance in data collection. Their expertise and suggestions have played a crucial role in the success of this project.


\bibliography{acl_latex}

\begin{thebibliography}{28}
\expandafter\ifx\csname natexlab\endcsname\relax\def\natexlab#1{#1}\fi

\bibitem[{Anil et~al.(2023)Anil, Dai, Firat, Johnson et~al.}]{anil2023palm}
Rohan Anil, Andrew~M Dai, Orhan Firat, Melvin Johnson, et~al. 2023.
\newblock Palm 2 technical report.
\newblock \emph{arXiv preprint arXiv:2305.10403}.

\bibitem[{Baldazzi et~al.(2023)Baldazzi, Bellomarini, Ceri, Colombo, Gentili, and Sallinger}]{baldazzi2023finetuning}
Teodoro Baldazzi, Luigi Bellomarini, Stefano Ceri, Andrea Colombo, Andrea Gentili, and Emanuel Sallinger. 2023.
\newblock Fine-tuning large enterprise language models via ontological reasoning.
\newblock \emph{arXiv preprint arXiv:2306.10723}.

\bibitem[{Brown et~al.(2020)Brown, Mann, Ryder, Subbiah et~al.}]{brown2020gpt3}
Tom Brown, Benjamin Mann, Nick Ryder, Melanie Subbiah, et~al. 2020.
\newblock Language models are few-shot learners.
\newblock \emph{Advances in neural information processing systems}, 33:1877--1901.

\bibitem[{Houlsby et~al.(2019)Houlsby, Giurgiu, Jastrzebski, Morrone et~al.}]{houlsby2019parameter}
Neil Houlsby, Andrei Giurgiu, Stanislaw Jastrzebski, Bruna Morrone, et~al. 2019.
\newblock Parameter-efficient transfer learning for nlp.
\newblock In \emph{International Conference on Machine Learning}, pages 2790--2799.

\bibitem[{Hu et~al.(2021)Hu, Shen, Wallis, Allen-Zhu et~al.}]{hu2021lora}
Edward~J Hu, Yelong Shen, Phillip Wallis, Zeyuan Allen-Zhu, et~al. 2021.
\newblock Lora: Low-rank adaptation of large language models.
\newblock \emph{arXiv preprint arXiv:2106.09685}.

\bibitem[{Imani et~al.(2023)Imani, Du, and Shrivastava}]{imani2023mathprompter}
Shima Imani, Liang Du, and Harsh Shrivastava. 2023.
\newblock Mathprompter: Mathematical reasoning using large language models.
\newblock \emph{arXiv preprint arXiv:2303.05398}.

\bibitem[{Li et~al.(2023)Li, Mellou, Zhang, Pathuri, and Menache}]{li2023supplychain}
Beibin Li, Konstantina Mellou, Bo~Zhang, Jeevan Pathuri, and Ishai Menache. 2023.
\newblock Large language models for supply chain optimization.
\newblock \emph{arXiv preprint arXiv:2307.03875}.

\bibitem[{Liu et~al.(2023)Liu, Yuan, Fu, Jiang et~al.}]{liu2023prompt}
Pengfei Liu, Weizhe Yuan, Jinlan Fu, Zhengbao Jiang, et~al. 2023.
\newblock Pre-train, prompt, and predict: A systematic survey of prompting methods in natural language processing.
\newblock \emph{ACM Computing Surveys}, 55(9):1--35.

\bibitem[{Madaan et~al.(2023)Madaan, Tandon, Gupta, Hallinan et~al.}]{madaan2023selfrefine}
Aman Madaan, Niket Tandon, Prakhar Gupta, Skyler Hallinan, et~al. 2023.
\newblock Self-refine: Iterative refinement with self-feedback.
\newblock \emph{arXiv preprint arXiv:2303.17651}.

\bibitem[{OpenAI(2023)}]{openai2023gpt4}
OpenAI. 2023.
\newblock Gpt-4 technical report.
\newblock \emph{arXiv preprint arXiv:2303.08774}.

\bibitem[{Ouyang et~al.(2022)Ouyang, Wu, Jiang, Almeida et~al.}]{ouyang2022training}
Long Ouyang, Jeffrey Wu, Xu~Jiang, Diogo Almeida, et~al. 2022.
\newblock Training language models to follow instructions with human feedback.
\newblock \emph{Advances in Neural Information Processing Systems}, 35:27730--27744.

\bibitem[{Patel et~al.(2021)Patel, Bhattamishra, and Goyal}]{patel2021mwp}
Arkil Patel, Satwik Bhattamishra, and Navin Goyal. 2021.
\newblock Are nlp models really able to solve simple math word problems?
\newblock \emph{arXiv preprint arXiv:2103.07191}.

\bibitem[{Qiao et~al.(2022)Qiao, Ou, Zhang, Chen et~al.}]{qiao2022reasoning}
Shuofei Qiao, Yixin Ou, Ningyu Zhang, Xiang Chen, et~al. 2022.
\newblock Reasoning with language model prompting: A survey.
\newblock \emph{arXiv preprint arXiv:2212.09597}.

\bibitem[{Ramamonjison et~al.(2022)Ramamonjison, Li, Yu, He et~al.}]{ramamonjison2022optimization}
Rindranirina Ramamonjison, Haley Li, Timothy~T Yu, Shiqi He, et~al. 2022.
\newblock Augmenting operations research with auto-formulation of optimization models from problem descriptions.
\newblock \emph{arXiv preprint arXiv:2209.15565}.

\bibitem[{Ramamonjison et~al.(2023)Ramamonjison, Yu, Li, Li et~al.}]{ramamonjison2023nl4opt}
Rindranirina Ramamonjison, Timothy~T Yu, Raymond Li, Haley Li, et~al. 2023.
\newblock Nl4opt competition: Formulating optimization problems based on their natural language descriptions.
\newblock \emph{arXiv preprint arXiv:2303.08233}.

\bibitem[{Shen et~al.(2023)Shen, Song, Tan, Li et~al.}]{shen2023hugginggpt}
Yongliang Shen, Kaitao Song, Xu~Tan, Dongsheng Li, et~al. 2023.
\newblock Hugginggpt: Solving ai tasks with chatgpt and its friends in huggingface.
\newblock \emph{arXiv preprint arXiv:2303.17580}.

\bibitem[{Silver et~al.(2017)Silver, Schrittwieser, Simonyan, Antonoglou et~al.}]{silver2017competitivestrategy}
David Silver, Julian Schrittwieser, Karen Simonyan, Ioannis Antonoglou, et~al. 2017.
\newblock Mastering the game of go without human knowledge.
\newblock \emph{nature}, 550(7676):354--359.

\bibitem[{Sun et~al.(2022)Sun, Shao, Qian, Huang, and Qiu}]{sun2022black}
Tianxiang Sun, Yunfan Shao, Hong Qian, Xuanjing Huang, and Xipeng Qiu. 2022.
\newblock Black-box tuning for language-model-as-a-service.
\newblock In \emph{International Conference on Machine Learning}, pages 20841--20855.

\bibitem[{Touvron et~al.(2023)Touvron, Lavril, Izacard, Martinet et~al.}]{touvron2023llama}
Hugo Touvron, Thibaut Lavril, Gautier Izacard, Xavier Martinet, et~al. 2023.
\newblock Llama: Open and efficient foundation language models.
\newblock \emph{arXiv preprint arXiv:2302.13971}.

\bibitem[{Vaswani et~al.(2017)Vaswani, Shazeer, Parmar, Uszkoreit et~al.}]{NIPS2017transformer}
Ashish Vaswani, Noam Shazeer, Niki Parmar, Jakob Uszkoreit, et~al. 2017.
\newblock Attention is all you need.
\newblock \emph{Advances in neural information processing systems}, 30.

\bibitem[{Wang et~al.(2022)Wang, Wei, Schuurmans, Le et~al.}]{wang2022selfconsistency}
Xuezhi Wang, Jason Wei, Dale Schuurmans, Quoc Le, et~al. 2022.
\newblock Self-consistency improves chain of thought reasoning in language models.
\newblock \emph{arXiv preprint arXiv:2203.11171}.

\bibitem[{Wei et~al.(2022)Wei, Wang, Schuurmans, Bosma et~al.}]{wei2022chain}
Jason Wei, Xuezhi Wang, Dale Schuurmans, Maarten Bosma, et~al. 2022.
\newblock Chain-of-thought prompting elicits reasoning in large language models.
\newblock \emph{Advances in Neural Information Processing Systems}, 35:24824--24837.

\bibitem[{Xie et~al.(2020)Xie, Dai, Chen, Dai et~al.}]{xie2020logistictransport}
Yujia Xie, Hanjun Dai, Minshuo Chen, Bo~Dai, et~al. 2020.
\newblock Differentiable top-k with optimal transport.
\newblock \emph{Advances in Neural Information Processing Systems}, 33:20520--20531.

\bibitem[{Yang et~al.(2023)Yang, Swope, Gu, Chalamala et~al.}]{yang2023leandojo}
Kaiyu Yang, Aidan~M Swope, Alex Gu, Rahul Chalamala, et~al. 2023.
\newblock Leandojo: Theorem proving with retrieval-augmented language models.
\newblock \emph{arXiv preprint arXiv:2306.15626}.

\bibitem[{Yao et~al.(2022)Yao, Zhao, Yu, Du et~al.}]{yao2022react}
Shunyu Yao, Jeffrey Zhao, Dian Yu, Nan Du, et~al. 2022.
\newblock React: Synergizing reasoning and acting in language models.
\newblock \emph{arXiv preprint arXiv:2210.03629}.

\bibitem[{Yao et~al.(2023)Yao, Aminabadi, Ruwase, Rajbhandari et~al.}]{yao2023deepspeed}
Zhewei Yao, Reza~Yazdani Aminabadi, Olatunji Ruwase, Samyam Rajbhandari, et~al. 2023.
\newblock Deepspeed-chat: Easy, fast and affordable rlhf training of chatgpt-like models at all scales.
\newblock \emph{arXiv preprint arXiv:2308.01320}.

\bibitem[{Ye et~al.(2020)Ye, Pei, Wang, Chen et~al.}]{ye2020financialinvestment}
Yunan Ye, Hengzhi Pei, Boxin Wang, Pin-Yu Chen, et~al. 2020.
\newblock Reinforcement-learning based portfolio management with augmented asset movement prediction states.
\newblock In \emph{Proceedings of the AAAI Conference on Artificial Intelligence}, volume~34, pages 1112--1119.

\bibitem[{Zhou et~al.(2023)Zhou, Li, and Liu}]{zhou2023llm}
Xuanhe Zhou, Guoliang Li, and Zhiyuan Liu. 2023.
\newblock Llm as dba.
\newblock \emph{arXiv preprint arXiv:2308.05481}.

\end{thebibliography}





\end{document}